\documentclass[12pt]{amsart}
\usepackage{amsmath,amscd,amssymb,amsfonts,graphics,mathrsfs}
\usepackage{hyperref}
\newcommand{\hl}{\hyperlink}
\newcommand{\htt}{\hypertarget}
\newcommand{\h}{\hbox}
\newcommand{\q}{\quad}

\newcommand{\bs}{\par\bigskip}
\newcommand{\ms}{\par\medskip}
\newcommand{\sk}{\par\smallskip}
\newcommand{\bsn}{\par\bigskip\noindent}
\newcommand{\msn}{\par\medskip\noindent}

\newcommand{\ges}{\geqslant}
\newcommand{\les}{\leqslant}
\newcommand{\1}{\hskip1pt}
\newcommand{\mcap}{\hbox{$\bigcap$}}

\newcommand{\msum}{\hbox{$\sum$}}
\newcommand{\mopl}{\hbox{$\bigoplus$}}
\newcommand{\mprod}{\hbox{$\prod$}}
\newcommand{\B}{{\mathcal B}}
\newcommand{\D}{{\mathcal D}}
\newcommand{\E}{{\mathcal E}}
\newcommand{\F}{{\mathscr F}}
\newcommand{\G}{{\mathcal G}}
\newcommand{\Hc}{{\mathcal H}}
\newcommand{\Lc}{{\mathcal L}}
\newcommand{\Nc}{{\mathcal N}}
\newcommand{\M}{{\mathcal M}}
\newcommand{\OO}{{\mathcal O}}

\newcommand{\Hs}{{\mathscr H}}
\newcommand{\Ls}{{\mathscr L}}
\newcommand{\Q}{{\mathbb Q}}
\newcommand{\C}{{\mathbb C}}
\newcommand{\N}{{\mathbb N}}

\newcommand{\DD}{{\mathbf D}}
\newcommand{\Z}{{\mathbb Z}}
\newcommand{\ee}{{\mathbf e}}
\newcommand{\mm}{{\mathfrak m}}
\newcommand{\de}{\delta}
\newcommand{\alt}{\widetilde{\alpha}}
\newcommand{\Bt}{\widetilde{\mathcal B}}
\newcommand{\Ht}{\widetilde{H}}
\newcommand{\Vt}{\widetilde{V}}
\newcommand{\Hst}{\widetilde{\mathscr H}}
\newcommand{\al}{\alpha}
\newcommand{\la}{\lambda}
\newcommand{\eBS}{e^{\rm BS}}
\newcommand{\Gr}{{\rm Gr}}
\newcommand{\dd}{\partial}
\newcommand{\bl}{\bigl}
\newcommand{\br}{\bigr}
\newcommand{\pl}{\1{+}\1}
\newcommand{\mi}{\1{-}\1}
\newcommand{\eq}{\,{=}\,}
\newcommand{\less}{\,{\leqslant}\,}
\newcommand{\gess}{\,{\geqslant}\,}

\newcommand{\ssb}{\raise.15ex\h{${\scriptscriptstyle\bullet}$}}
\newcommand{\ssc}{\,\raise.15ex\h{${\scriptstyle\circ}$}\,}
\newcommand{\onto}{\twoheadrightarrow}
\newcommand{\into}{\hookrightarrow}
\newcommand{\simto}{\,\,\rlap{\hskip1.5mm\raise1.4mm\hbox{$\sim$}}\hbox{$\longrightarrow$}\,\,}
\begin{document}
\title[Brian\c con-Skoda exponents]{Brian\c con-Skoda exponents and the maximal root of reduced Bernstein-Sato polynomials}
\author[S.-J. Jung]{Seung-Jo Jung}
\address{S.-J. Jung : Department of Mathematics Education, and Institute of Pure and Applied Mathematics, Jeonbuk National University, Jeonju, 54896, Korea}
\email{seungjo@jbnu.ac.kr}
\author[I.-K. Kim]{In-Kyun Kim}
\address{I.-K. Kim : Department of Mathematics, Yonsei University, Seoul, Korea}
\email{soulcraw@gmail.com}
\author[M. Saito]{Morihiko Saito}
\address{M. Saito : RIMS Kyoto University, Kyoto 606-8502 Japan}
\email{msaito@kurims.kyoto-u.ac.jp}
\author[Y. Yoon]{Youngho Yoon}
\address{Y. Yoon : Department of Mathematics, Chungnam National University, 99 Daehak-ro, Daejeon 34134, Korea}
\email{mathyyoon@skku.edu}
\thanks{This work was partially supported by NRF grant funded by the Korea government(MSIT) (the first author: NRF-2021R1C1C1004097, the second author: NRF-2020R1A2C4002510, and the fourth author: NRF-2020R1C1C1A01006782). The third author was partially supported by JSPS Kakenhi 15K04816.}
\begin{abstract} For a holomorphic function $f$ on a complex manifold $X$, the Brian\c con-Skoda exponent $e^{\rm BS}(f)$ is the smallest integer $k$ with $f^k\in(\partial f)$ (replacing $X$ with a neighborhood of $f^{-1}(0)$), where $(\partial f)$ denotes the Jacobian ideal of $f$. It is shown that $e^{\rm BS}(f)\le d_X$ $(:=\dim X)$ by Brian\c con-Skoda. We prove that $e^{\rm BS}(f)\le[d_X-2\widetilde{\alpha}_f]+1$ with $-\widetilde{\alpha}_f$ the maximal root of the reduced Bernstein-Sato polynomial $b_f(s)/(s+1)$, assuming the latter exists (shrinking $X$ if necessary). This implies for instance that $e^{\rm BS}(f)\le d_X-2$ in the case $f^{-1}(0)$ has only rational singularities, that is, if $\widetilde{\alpha}_f>1$.
\end{abstract}
\maketitle
\centerline{\bf Introduction}
\bsn
Let $f$ be a holomorphic function on a complex manifold $X$ of dimension $d_X$. Set $Z:=f^{-1}(0)$. We assume the Bernstein-Sato polynomial $b_f(s)$ exists (shrinking $X$ if necessary). This coincides with the least common multiple of the local Bernstein-Sato polynomials $b_{f,x}(s)$ ($x\in Z$). Let $\alt_f$ be the maximal root of the reduced Bernstein-Sato polynomial $b_f(s)/(s{+}1)$ up to sign. This is a positive rational number assuming $Z$ is singular, see \cite{Ka}. It is also known (see \cite[Thm.\,0.4]{mic}) that in the singular case we have
\htt{1}{}
$$\alt_f\less d_X/2.
\leqno(1)$$
\sk
The {\it Brian\c con-Skoda exponent\1} is defined by
$$\eBS(f):=\min\bl\{k\in\N\mid f^k\in(\dd f)\br\},$$
replacing $X$ with an open neighborhood of $f^{-1}(0)$ if necessary (so that the $f^{-1}(t)$ are nonsingular for $t\ne 0$). Here $(\dd f)\subset\OO_X$ is the Jacobian ideal generated locally by the partial derivatives of $f$. By definition this invariant is essentially local; more precisely, we have
\htt{2}{}
$$\eBS(f)={\rm max}\{\eBS(f,x)\}_{x\in Z},
\leqno(2)$$
where $\eBS(f,x)$ is defined by using the ideal $(\dd f)_x\subset\OO_{X,x}$. We have $f^k\in(\dd f)_x$ for $k\,{\gg}\,0$, $x\in Z$ by an analytic version of Hilbert's Nullstellensatz. In the isolated singularity case, we have $\eBS(f)=1$ if and only if $f$ is a weighted homogeneous polynomial, see \cite{SaK}. So this number measures how far the function is away from weighted homogeneous polynomials.
\sk
In general it is proved that $\eBS(f)\les d_X$, that is, $f\1^{d_X}\in(\dd f)$ by Brian\c con-Skoda \cite{BrSk}. In this paper we show the following.
\par\htt{T1}{}\msn
{\bf Theorem~1.} {\it In the above notation, we have
\htt{3}{}
$$\eBS(f)\les[\1 d_X\mi 2\1\widetilde{\alpha}_f\1]\pl 1,
\leqno(3)$$
or equivalently,}
\htt{4}{}
$$f^k\in(\dd f)\q\h{\it if}\q k>d_X\mi 2\1\alt_f.
\leqno(4)$$
\sk
We prove Theorem~\hl{T1}{1} by using the {\it microlocal $V$-filtration\1} $\Vt$ on $\OO_X$, see \cite{mic}.
In the first proof, we calculated the {\it generating level\1} of each monodromy eigensheaf of the underlying filtered $\D$-modules of the vanishing cycle Hodge module, see \hl{2.1}{2.1} below.
(We use ``generating level" rather than ``generation level", since the latter may mean how much something is generated. The term ``level" comes from the ``Hodge level" of a Hodge structure (taking its direct image by a closed embedding of a point into a complex manifold).
For the proof of Theorem~\hl{T1}{1}, it is essential to determine the number $q(\M_{\la},F)$ in the notation of \hl{2.1}{2.1} below. It turns out that it is better to show the following.
\par\htt{T2}{}\msn
{\bf Theorem~2.} {\it We have the inclusion}
\htt{5}{}
$$\Vt^{\al}\OO_X\subset(\dd f)\q\h{if}\q\al>d_X\mi\alt_f.
\leqno(5)$$
\ms
This implies immediately Theorem~\hl{T1}{1} using the inclusions $f\1\Vt^{\al}\OO_X\subset\Vt^{\al+1}\OO_X$ ($\al\in\Q$), see (\hl{1.1.7}{1.1.7}) below. Note that we have the equality (see for instance \cite[(1.3.8)]{hi})
\htt{6}{}
$$\alt_f=\min\bl\{\al\in\Q\,\big|\,\Gr_{\Vt}^{\al}\OO_X\ne 0\br\}.
\leqno(6)$$
\sk
It has been remarked by M.~Musta\c t\u a that the proof of the assertion (\hl{4}{4}) in the {\it algebraic case\1} can be reduced to the isolated singularity case, see Remark\,\,\hl{R2.2a}{2.2a} below.
(It is not necessarily trivial to extend an argument of mixed Hodge modules in the algebraic case to the corresponding analytic case if the {\it derived categories of mixed Hodge modules\1} is used.)
\sk
In the isolated singularity case, Theorem~\hl{T1}{1} follows from mixed Hodge theory for vanishing cohomology related to Brieskorn lattices (see \cite{Bri}, \cite{Mal}, \cite{bl}, \cite{St}, \cite{SS}, \cite{Va2}, etc.), where $\alt_f$ coincides with the minimal spectral number, see also \cite[Rem.\,3.4b]{JKSY2} and 2.2 below. Theorem~\hl{T1}{1} may be useful for instance when one tries to calculate the sum of Milnor numbers of singular points of $Z$ using a computer program like Singular \cite{Sing} in the case $f$ is a polynomial of many variables provided that we know already that $\alt_f$ is large for any singular point of $Z$ (for instance, if $Z$ has only rational singularities, where $\alt_f>1$), see also \cite[Rem.\,6.2]{DiSt}, \cite[A.2]{wh}, \cite[Remark after Ex.\,5.2]{nwh}, etc.
\sk
In the isolated singularity case it is known that
\htt{7}{}
$${\rm NO}(f)\les\eBS(f),
\leqno(7)$$
where ${\rm NO}(f)$ is the {\it nilpotence order\1} of $N:=\log T_u$ with $T=T_sT_u$ the Jordan decomposition of the monodromy $T$ on the vanishing cohomology of $f$, that is, the largest size of Jordan blocks, see \cite{Sch1}, \cite{Sch2}, \cite{Va1}, \cite{Va2}, \cite{vSt}, and Remark\,\,\hl{R2.2b}{2.2b} below. (In the the non-isolated singularity case, the inequality (\hl{7}{7}) does not hold even if we define ${\rm NO}(f)$ using the vanishing cohomology groups $\Ht^{\ssb}(F_{\!f,x},\Q)$, where $F_{\!f,x}$ denotes the Milnor fiber of $f$ at $x\in Z$, see Remark\,\,\hl{R2.2d}{2.2d} below.) Combining (\hl{7}{7}) with (\hl{3}{3}) in Theorem~\hl{T1}{1}, we get the inequality
\htt{8}{}
$${\rm NO}(f)\les[\1 d_X\mi 2\1\widetilde{\alpha}_f\1]\pl 1,
\leqno(8)$$
which is more or less known related to the theory of spectral pairs \cite{St}. The inequalities (\hl{3}{3}) and (\hl{8}{8}) are optimal as is seen in the case of polynomials
\htt{9}{}
$$f_{a,b}=\mprod_{i=1}^{d_X}\,x_i^a+\msum_{i=1}^{d_X}\,x_i^b\q(a\gess 2,\,b\,{>}\,d_Xa,\,d_X\gess 2),
\leqno(9)$$
generalizing examples in \cite{MalL}, \cite{BrSk} (where $a\eq 2$ or 3, see also \cite{vSt}). Here we have $\alt_{f_{a,b}}\eq 1/a\les 1/2$ (see \cite{EL}, \cite{exp}), and ${\rm NO}(f)\eq d_X$, see \cite{Sta}, etc. (The last equality follows from the descent theorem of nearby cycle formula for Newton non-degenerate functions, see \cite[Cor.\,3]{des}.)
Combining these with the Thom-Sebastiani type theorem explained just below, it is possible to produce examples such that the equality holds in (\hl{8}{8}) and moreover ${\rm NO}(f)$ is equal to any given integer in $[1,d_X]$ for any dimension $d_X\gess2$, replacing $X,f_{a,b}$ with $X{\times}Y,f_{a,b}\pl g$, where $g:=\msum_{j=1}^{d_Y}\,y_j^2$.
\sk
As for the relation to the Thom-Sebastiani type theorem, set $h=f\pl g$ on $X{\times}Y$, where $g$ is a holomorphic function on a complex manifold $Y$ with $\eBS(g)=1$, that is, $g\in(\dd g)$. Then
\htt{10}{}
$$\eBS(h)=\eBS(f),\q\alt_h=\alt_f\pl\1\alt_g,
\leqno(10)$$
see for instance \cite[Thm.\,0.8]{mic} for the last equality. The inequality (\hl{3}{3}) for $f$ is stronger than that for $h$ by (\hl{1}{1}), and these are equivalent if and only if the equality holds in (\hl{1}{1}) for $g$, that is, $\alt_g\eq d_Y/2$, or equivalently, $g^{-1}(0)$ has only ordinary double points.
\sk
Theorem~\hl{T1}{1} also implies the following.
\par\htt{C1}{}\msn
{\bf Corollary\,\,1.} {\it If $Z$ has only rational singularities, we have}
\htt{11}{}
$$\eBS(f)\les d_X\mi 2.
\leqno(11)$$
\ms
Note that $Z$ has only rational singularities if and only if $\alt_f>1$, see \cite[Thm\,0.4]{rat}.
\sk
In Section 1 we review some basics of microlocal $V$-filtration and generating level.
In Section 2 we prove Theorem~\hl{T2}{2}, which implies Theorem~\hl{T1}{1} immediately, and give some remarks.
\sk
We thank M.~Musta\c t\u a for a useful comment about the reduction of the proof of Theorem~\hl{T1}{1} to the isolated singularity case when $X,f$ are algebraic.
We are grateful to the referee for careful reading of the paper.
\bs\bs
\vbox{\centerline{\bf 1. Preliminaries}
\bsn
In this section we review some basics of microlocal $V$-filtration and generating level.}
\msn
{\bf 1.1.~Microlocal $V$-filtration.} Let $f$ be a holomorphic function on a complex manifold $X$. We assume ${\rm Sing}\,f\subset f^{-1}(0)$ replacing $X$ with an open neighborhood of $f^{-1}(0)$ if necessary. Let $i_f:X\into X{\times}\C$ be the graph embedding with $t$ the coordinate of $\C$. Set
\htt{1.1.1}{}
$$(\B_f,F):=(i_f)_*^{\D}(\OO_X,F).
\leqno(1.1.1)$$
The right-hand side is the direct image as a filtered $\D$-module, see \cite{mhp}, \cite{ypg}, etc. It has the filtration $V$ along $t\eq 0$ indexed by $\Q$, see \cite{Ka2}, \cite{Mal2}, \cite{bl}, \cite{mhp}, etc. (Note that $V_{\al}=V^{-\al}$.)
\sk
Let $(\Bt_f,F)$ be the {\it algebraic partial microlocalization\1} of $(\B_f,F)$, see \cite{mic}, \cite{MSS}, etc. By definition $\B_f$ and $\Bt_f$ are {\it freely generated\1} over $\OO_X[\dd_t]$ and $\OO_X[\dd_t,\dd_t^{-1}]$ respectively by $\de(t{-}f)$ satisfying the relations
\htt{1.1.2}{}
$$t\de(t{-}f)=f\de(t{-}f),\q\dd_{x_i}\de(t{-}f)=-f_i\1\dd_t\1\de(t{-}f),
\leqno(1.1.2)$$
with $x_i$ local coordinates of $X$ and $f_i:=\dd_{x_i}f$.
\sk
We define the filtration $V$ on $\Bt_f$ as in \cite[2.1.3]{mic}\,:
\htt{1.1.3}{}
$$V^{\al}\Bt_f:=\begin{cases}V^{\al}\B_f+\OO_X[\dd_t^{-1}]\dd_t^{-1}\de(t{-}f)\q\q(\al\less 1),\\ \dd_t^{-j}V^{\al-j}\B_f\q\q\,\,\,\,\,(\al\,{>}\,1,\,\,\al\1{-}j\in(0,1\1]).
\end{cases}
\leqno(1.1.3)$$
\sk
In this paper the Hodge filtration $F$ is indexed as in the case of filtered right $\D$-modules, in particular,
\htt{1.1.4}{}
$$\Gr_p^F\Bt_f=\OO_X\,\dd_t^{\,p+d_X}\1\de(t{-}f)\q(p\in\Z).
\leqno(1.1.4)$$
\sk
We define the {\it microlocal $V\!$-filtration\1} on $\OO_X$ by
\htt{1.1.5}{}
$$(\OO_X,\Vt):=\Gr^F_{-d_X}(\Bt_f,V),
\leqno(1.1.5)$$
By the definition of $V$-filtration we have the inclusions
\htt{1.1.6}{}
$$t\1V^{\al}\Bt_f\subset V^{\al+1}\Bt_f,\q \dd_{x_i}V^{\al}\Bt_f\subset V^{\al}\Bt_f\q(\al\in\Q).
\leqno(1.1.6)$$
Using (\hl{1.1.2}{1.1.2}), these imply that
\htt{1.1.7}{}
$$f\1\Vt^{\al}\OO_X\subset\Vt^{\al+1}\OO_X,\q f_i\1\Vt^{\al}\OO_X\subset\Vt^{\al+1}\OO_X\,\,\,(i\in[1,d_X]).
\leqno(1.1.7)$$
\sk
Let $(\M_{\la},F)$ be the $\la$-eigensheaf of the monodromy on the underlying filtered left $\D_X$-module of the vanishing cycle Hodge module (see \cite{mhp}, \cite{mhm}) which is denoted by $\varphi_f\Q_{h,X}[d_Z]$ in this paper. By \cite[(5.1.3.3)]{mhp}, \cite[(2.1.4), (2.2.1)]{mic} (see also Remark\,\,\hl{R1.1a}{1.1a} below) we have the following.
\par\htt{P1.1}{}\msn
{\bf Proposition\,\,1.1.} {\it There are canonical isomorphisms
\htt{1.1.8}{}
$$\Gr_V^{\alpha}(\Bt_f,F)=(\M_{\la},F)\q\bl(\al\in(-1,0],\,\,\la=e^{-2\pi i\al}\br),
\leqno(1.1.8)$$
\vskip-6mm
\htt{1.1.9}{}
$$\hskip-4mm\dd_t^j:(\Bt_f;F,V)\simto(\Bt_f;F[-j],V[-j])\q(j\in\Z).
\leqno(1.1.9)$$
with $(F[-j])_p:=F_{p+j}$, $(V[-j])^{\al}:=V^{\al-j}\,\,\,($in a compatible way with $F^p=F_{-p})$.}
\ms
Using the canonical isomorphisms
\htt{1.1.10}{}
$$\Gr_V^{\al}\Gr^F_p\Bt_f=\Gr^F_p\Gr_V^{\al}\Bt_f\q(p\in\Z,\,\,\al\in\Q),
\leqno(1.1.10)$$
together with the last equality in (\hl{1.1.2}{1.1.2}), we then get the following.
\par\htt{C1.1a}{}\msn
{\bf Corollary\,\,1.1a.} {\it There are canonical isomorphisms for $\al\in\Q$}\,:
\htt{1.1.11}{}
$$\Gr_{\Vt}^{\al}\OO_X=\Gr^F_{-p}\M_{\la}\q\bl([d_X{-}\al]=p,\,\,\la=e^{-2\pi i\al}\br),
\leqno(1.1.11)$$
\par\htt{C1.1b}{}\msn
{\bf Corollary\,\,1.1b.} {\it The graded action of $-f_i$ on $\Gr_{\Vt}^{\ssb}\OO_X$ is identified with the action of $\Gr^F_1\dd_{x_i}$ on $\Gr^F_{\ssb}\M_{\la}$.}
\par\htt{R1.1a}{}\msn
{\bf Remark\,\,1.1a.} By \cite[(3.2.1.3), (5.1.3.3)]{mhp} we have
\htt{1.1.12}{}
$$\aligned\psi_f(\B_f,F)&=\mopl_{0<\al\les1}\,\Gr_V^{\al}(\B_f,F[1]),\\ \varphi_f(\B_f,F)&=\mopl_{-1<\al\les0}\,\Gr_V^{\al}(\B_f,F),\endaligned
\leqno(1.1.12)$$
where the action of $\dd_tt{-}\al$ corresponds to $-N/2\pi i$ on $\psi_f\C_X[d_X{-}1]$, $\varphi_f\C_X[d_X{-}1]$ by the de Rham functor, see \cite{Mal2}, etc. Note that the non-unipotent monodromy part of the nearby and vanishing cycles $\psi_{f,\ne1}(\B_f,F)$ and $\varphi_{f,\ne1}(\B_f,F)$ coincide by definition.
\ms
For the proof of the main theorems we need also the following.
\par\htt{L1.1}{}\msn
{\bf Lemma\,\,1.1.} {\it For any $\al\in\Q$, the graded $\Gr^F_{\ssb}\D_X$-modules $\Gr^F_{\ssb}V^{\al}\B_f$ are locally finitely generated.}
\msn
{\it Proof.} Since mixed Hodge modules are stable under the nearby and vanishing cycle functors (which are defined by using the graded pieces of the filtration $V$ for the underlying filtered $\D$-modules as in (\hl{1.1.12}{1.1.12})), the graded $\Gr^F_{\ssb}\D_X$-modules $\Gr^F_{\ssb}\Gr_V^{\al}\B_f$ are locally finitely generated for any $\al\in\Q$, see \cite[(5.1.3.3), (5.1.6.2)]{mhp}. It is then enough to consider the case $\al>0$ using the following short exact sequences inductively for the case $\al\less 0$\,:
\htt{1.1.13}{}
$$0\to(V^{>\al}\B_f,F)\to(V^{\al}\B_f,F)\to(\Gr_V^{\al}\B_f,F)\to0.
\leqno(1.1.13)$$
\sk
By \cite[(3.2.1.2)]{mhp}, we have the isomorphisms
\htt{1.1.14}{}
$$f\1\cdot:\Gr^F_pV^{\al}\B_f\simto\Gr^F_pV^{\al+1}\B_f\q(\al\,{>}\,0,\,p\in\Z),
\leqno(1.1.14)$$
since the action of $t$ on $\Gr^F_p\B_f$ is given by multiplication by $f$, which is denoted by $f\1\cdot$ in (\hl{1.1.14}{1.1.14}). (Note that (\hl{1.1.14}{1.1.14}) does not hold with $\B_f$ replaced by $\Bt_f$.)
\sk
We then see that the graded $\Gr^F_{\ssb}\D_X$-modules
$$\Gr^F_{\ssb}V^{\al}\B_f/f(\Gr^F_{\ssb}V^{\al}\B_f)=\Gr^F_{\ssb}(V^{\al}\B_f/V^{\al+1}\B_f)$$
are locally finitely generated over $\Gr^F_{\ssb}\D_X\eq\OO_X[\xi_1,\dots,\xi_n]$ with $\xi_i:=\Gr^F_1\dd_{x_i}$ and $n\eq d_X$. This implies that the graded $\OO_X$-modules
$$\G^{\al}_{\ssb}:=\Gr^F_{\ssb}V^{\al}\B_f\otimes_{\Gr^F_{\ssb}\D_X}(\Gr^F_{\ssb}\D_X/\Gr^F_{\ssb>0}\D_X)$$
with $\Gr^F_{\ssb}\D_X/\Gr^F_{\ssb>0}\D_X\eq\OO_X$ are locally finitely generated over $\OO_X$ for each $\al\,{>}\,0$ (since so are their tensor products with $\OO_X/(f)$). Indeed, the $\G^{\al}_p$ are coherent for any $p\in\Z$, and are supported on $f^{-1}(0)$ except for $p\eq{-}d_X$ (since at any point of $X\setminus f^{-1}(0)$, some $f_i$ is invertible). Using Nakayama's lemma (see for instance \cite[Cor.\,4.8]{Ei1}), we then conclude that $\G^{\al}_p\eq0$ except for a finite number of $p$ locally on $X$ for each $\al\,{>}\,0$ (since its tensor product with $\OO_X/(f)$ vanishes).
\sk
The last assertion implies the following equalities for $p\,{\gg}\,0$ locally on $X$ (depending on $\al\,{>}\,0$)
\htt{1.1.15}{}
$$\Gr^F_pV^{\al}\B_f=(\Gr^F_1\D_X)\1\Gr^F_{p-1}V^{\al}\B_f,
\leqno(1.1.15)$$
since $\Gr^F_q\D_X\eq(\Gr^F_1\D_X)\1\Gr^F_{q-1}\D_X$ for $q\,{>}\,0$. (This may be viewed as a graded version of Nakayama's lemma for graded $\Gr^F_{\ssb}\D_X$-modules. Here it seems enough to assume that a graded module is {\it lower bounded,} that is, its degree $p$ part vanishes if $p\,{\ll}\,0$. Finite generation condition as in \cite[Lem.\,1.4]{Ei2}) seems unnecessary assuming lower boundedness.) Since the $\Gr^F_pV^{\al}\B_f$ are coherent $\OO_X$-modules, the assertion then follows. This finishes the proof of Lemma\,\,\hl{L1.1}{1.1}.
\par\htt{C1.1c}{}\msn
{\bf Corollary\,\,1.1c.} {\it We have locally on $X$}
\htt{1.1.16}{}
$$\Vt^{\al}\OO_X=\msum_{i=1}^{d_X}\,f_i\1\Vt^{\al-1}\OO_X\q\h{if}\q\al\gg 0.
\leqno(1.1.16)$$
\msn
{\it Proof.} By (\hl{1.1.3}{1.1.3}--\hl{1.1.4}{4}) we have the isomorphisms for $\al\less 1$, $\,p\gess{-}d_X$\,:
\htt{1.1.17}{}
$$\Gr^F_pV^{\al}\Bt_f=\Gr^F_pV^{\al}\B_f.
\leqno(1.1.17)$$
So the assertion follows from Lemma\,\,\hl{L1.1}{1.1} using (\hl{1.1.9}{1.1.9}).
\ms
Using the $\Gr^F_{\ssb}$ of the short exact sequences (\hl{1.1.13}{1.1.13}), we get the following.
\par\htt{C1.1d}{}\msn
{\bf Corollary\,\,1.1d.} {\it For $\al_0\in\Q$, the equality {\rm(\hl{1.1.16}{1.1.16})} holds for any $\al\,{>}\,\al_0$ if and only if}
\htt{1.1.18}{}
$$\Gr_{\Vt}^{\al}\OO_X=\msum_{i=1}^{d_X}\,(\Gr_{\Vt}f_i)\1\Gr_{\Vt}^{\al-1}\OO_X\q\h{\it for any}\,\,\,\al\,{>}\,\al_0.
\leqno(1.1.18)$$
\msn
{\bf 1.2.~Generating level.} Let $(\M,F)$ be a holonomic filtered right $\D_X$-module. Put
\htt{1.2.1}{}
$$\aligned p(\M,F)&:=\min\bl\{p\in\Z\,\big|\,F_p\M\ne 0\br\},\\ q(\M,F)&:=\max\bl\{q\in\Z\,\big|\,(\Gr^F_{q-1}\M)\Gr^F_1\D_X\ne\Gr^F_q\M\br\},\\ \ell_{\rm gen}(\M,F)&:=q(\M,F)\mi p(\M,F).\endaligned
\leqno(1.2.1)$$
The last invariant is called the {\it generating level\1} of $\M$, see \cite[1.1]{mos}.
If $(\M,F)$ is the direct image of a Hodge structure by a closed embedding of a point, then the generating level coincides with the ``Hodge level" of the Hodge structure.
\sk
Assume there is locally a filtered free resolution
\htt{1.2.2}{}
$$(\Lc^{\ssb},F)\simto(\M,F),
\leqno(1.2.2)$$
where
\htt{1.2.3}{}
$$(\Lc^j,F)=(\E^j,F){\otimes}_{\OO_X}(\D_X,F),
\leqno(1.2.3)$$
with $\E^j$ coherent $\OO_X$-modules endowed with a finite filtration $F$ whose graded pieces are free $\OO_X$-modules of finite rank. We assume
\htt{1.2.4}{}
$$\Gr^F_p\1\E^j\eq0\q\h{if}\,\,\,j\ll 0\,\,\,\h{or}\,\,\,j>0.
\leqno(1.2.4)$$
Set
\htt{1.2.5}{}
$$\aligned p(\Lc^{\ssb},F)&:=\min\bl\{p\in\Z\,\big|\,\Gr^F_p\1\E^j\,{\ne}\,0\,\,\,(\exists\,j\in\Z)\br\},\\ q(\Lc^{\ssb},F)&:=\max\bl\{q\in\Z\,\big|\,\Gr^F_q\1\E^j\,{\ne}\,0\,\,\,(\exists\,j\in\Z)\br\},\\ \ell_{\rm res}(\Lc^{\ssb},F)&:=q(\Lc^{\ssb},F)\mi p(\Lc^{\ssb},F),\\ \ell_{\rm res}(\M,F)&:=\min\bl\{\ell_{\rm res}(\Lc^{\ssb},F)\,\big|\,(\Lc^{\ssb},F)\simto(\M,F)\br\}.\endaligned
\leqno(1.2.5)$$
Here $(\Lc^{\ssb},F)\simto(\M,F)$ runs over filtered free resolutions satisfying (\hl{1.2.3}{1.2.3}--\hl{1.2.4}{4}).
In this paper, $\ell_{\rm res}(\M,F)$ is called the {\it resolution level.} This is defined only locally (replacing $X$ with a sufficiently small neighborhood of each $x\in X$). Note that this definition is different from \cite[1.1]{mos}. (If $(\M,F)$ underlies a variation of Hodge structure on $X$, the resolution level coincides with the sum of the stalkwise Hodge level of $(\M,F)$ and $d_X$, see also Remark~\hl{R1.2}{1.2} below.)
\sk
By Lemma\,\,\hl{L1.2b}{1.2b} below we may assume
\htt{1.2.6}{}
$$p(\Lc^{\ssb},F)=p(\M,F).
\leqno(1.2.6)$$
\sk
We have the following.
\par\htt{L1.2a}{}\msn
{\bf Lemma\,\,1.2a.} {\it In the above notation and assumptions, the following inequalities hold\,}:
\htt{1.2.7}{}
$$p(\M,F)\ges p(\Lc^{\ssb},F),\q q(\M,F)\les q(\Lc^{\ssb},F),
\leqno(1.2.7)$$
\vskip-7mm
\htt{1.2.8}{}
$$\ell_{\rm gen}(\M,F)\les\ell_{\rm res}(\M,F).
\leqno(1.2.8)$$
\msn
{\it Proof.} It is enough to show the two inequalities in (\hl{1.2.7}{1.2.7}). The first inequality follows from (\hl{1.2.3}{1.2.3}--\hl{1.2.4}{4}). For the second one, note that $q(\M,F)$ coincides with
$$\min\bl\{q\in\Z\,\big|\,(\Gr^F_k\M)\Gr^F_1\D_X\eq\Gr^F_{k+1}\M\,\,(\forall\,k\gess q)\br\}.$$
The assertion is then proved using the commutative diagram
\htt{1.2.9}{}
$$\begin{array}{ccc}\Gr^F_k\Lc^0\otimes_{\OO_X}\Gr^F_1\D_X&\onto&\Gr^F_{k+1}\Lc^0\\ \rotatebox{270}{$\!\!\!\onto$\,\,\,}&&\rotatebox{270}{$\!\!\!\onto$\,\,\,}\\\Gr^F_k\M\otimes_{\OO_X}\Gr^F_1\D_X&\to&\Gr^F_{k+1}\M\end{array}
\leqno(1.2.9)$$
where the vertical and top surjectivity for $k\gess q(\Lc^{\ssb},F)$ follows respectively from (\hl{1.2.4}{1.2.4}) and (\hl{1.2.3}{1.2.3}). So Lemma\,\,\hl{L1.2a}{1.2a} follows.
\par\htt{L1.2b}{}\msn
{\bf Lemma\,\,1.2b.} {\it If there is locally a filtered free resolution $(1.2.2)$ satisfying $(1.2.3$--$4)$, then there is locally a filtered free resolution
\htt{1.2.10}{}
$$(\Lc'{}^{\ssb},F)\simto(\M,F),
\leqno(1.2.10)$$
satisfying $(1.2.3$--$4)$ and}
\htt{1.2.11}{}
$$p(\Lc'{}^{\ssb},F)=p(\M,F),\q q(\Lc'{}^{\ssb},F)=q(\Lc^{\ssb},F).
\leqno(1.2.11)$$
\msn
{\it Proof.} We have a finite filtration $G$ on $(\Lc^j,F)$ induced from the filtration $F$ on $\E^j$, that is,
\htt{1.2.12}{}
$$G_k\Lc^j=F_k\E^j\otimes_{\OO_X}\D_X,
\leqno(1.2.12)$$
see also \cite[2.1.10]{mhp}. We see that $\Gr^G_p(\Lc^{\ssb},F)$ is filtered acyclic if $p=p(\Lc^{\ssb},F)$ and $\Gr^F_p\Lc^{\ssb}$ is acyclic. Here the last condition is equivalent to that $p(\Lc^{\ssb},F)<p(\M,F)$. So the assertion follows by induction on the difference between $p(\M,F)$ and $p(\Lc^{\ssb},F)$ replacing $(\Lc^{\ssb},F)$ by $(\Lc^{\ssb},F)/(F_p\Lc^{\ssb},F)$ with $p=p(\Lc^{\ssb},F)$ in the case $p(\Lc^{\ssb},F)<p(\M,F)$. This finishes the proof of Lemma\,\,\hl{L1.2b}{1.2b}.
\par\htt{R1.2}{}\msn
{\bf Remark\,\,1.2.} The equality $q(\Lc'{}^{\ssb},F)=q(\M,F)$ does not necessarily hold in general; for instance, if $\M=\Omega^{d_X}_X$, which is locally isomorphic to $\D_X/\msum_i\,\dd_{x_i}\D_X$.
\bs\bs
\vbox{\centerline{\bf 2. Proof of the main theorems}
\bsn
In this section we prove Theorem~\hl{T2}{2}, which implies Theorem~\hl{T1}{1} immediately, and give some remarks.}
\par\htt{2.1}{}\msn
{\bf 2.1.~Proof of Theorem~\hl{T2}{2}.} For $\al\in\Q$, set
$$\ee(\al):=e^{2\pi i\al},\q\{\al\}:=\al\mi[\al],\q\{\al\}^{\!\vee}:=-\{-\al\},$$
that is, $\{\al\}^{\!\vee}\in(-1,0\1]$, $\al\mi\{\al\}^{\!\vee}\in\Z$ ($\al\in\Q$).
\sk
By Corollary\,\,\hl{C1.1a}{1.1a}, we have for $\al\in(-1,0\1]$\,:
\htt{2.1.1}{}
$$p(\M_{\ee(-\al)},F)=\begin{cases}-[d_X\mi\alt_f]&(\al\ges\{\alt_f\}^{\!\vee}),\\-[d_X\mi\alt_f]\pl 1&(\al<\{\alt_f\}^{\!\vee}).\end{cases}
\leqno(2.1.1)$$
\sk
It is known that the $\Gr^F_{\ssb}\M_{\la,x}$ are Cohen-Macaulay $\Gr^F_{\ssb}\D_{X,x}$-modules for $x\in X$, see \cite[Lemma\,\,5.1.13]{mhp}. Since the $\M_{\la}$ are holonomic $\D_X$-modules, this implies that locally on $X$, there are filtered free resolutions of length $d_X$\,:
\htt{2.1.2}{}
$$(\Nc_{\la}^{\ssb},F)\simto(\M_{\la},F),
\leqno(2.1.2)$$
using the Auslander-Buchsbaum formula for the projective dimension of $\Gr^F_{\ssb}\M_{\la,x}$ (see \cite{BrHe}, \cite{Ei1}) together with the graded Nakayama lemma. Here we may assume by Lemma\,\,\hl{L1.2b}{1.2b}
\htt{2.1.3}{}
$$p(\Nc_{\la}^{\ssb},F)=p(\M_{\la},F).
\leqno(2.1.3)$$
\sk
We have the self-duality isomorphisms (see \cite{dual})\,:
\htt{2.1.4}{}
$$\DD(\M_{\la},F)=\begin{cases}(\M_{\la^{-1}},F[d_X{-}1])&(\la\ne 1),\\ (\M_{\1 1},F[d_X])&(\la=1),\end{cases}
\leqno(2.1.4)$$
Here the $\DD(\M_{\la},F)$ can be calculated by using the free resolution (\hl{2.1.2}{2.1.2}). So we have a filtered free resolution
\htt{2.1.5}{}
$$(\Nc_{\la}^{\prime\ssb},F)\simto\DD(\M_{\la},F),
\leqno(2.1.5)$$
by putting first $\,(\Nc_{\la}^{\prime\ssb},F):=\DD(\Nc_{\la}^{\ssb},F)$ with $\DD(\Nc_{\la}^{\ssb},F)$ defined by
\htt{2.1.6}{}
$$\Hc om_{\OO_X}\bl((\Nc_{\la}^{\ssb},F),(\omega_X,F)\otimes_{\OO_X}(\D_X,F)\br)[d_X],
\leqno(2.1.6)$$
where $\Gr^F_p\omega_X=0$ ($p\ne 0$), see \cite{mhp}, \cite{ypg}, \cite{JKSY2}.
\sk
By (\hl{2.1.1}{2.1.1}), (\hl{2.1.3}{2.1.3}), (\hl{2.1.6}{2.1.6}), we have for $\al\in(-1,0\1]$
\htt{2.1.7}{}
$$q(\Nc_{\ee(-\al)}^{\prime\ssb},F)=\begin{cases}[d_X\mi\alt_f]&(\al\ges\{\alt_f\}^{\!\vee}),\\ [d_X\mi\alt_f]\mi 1&(\al<\{\alt_f\}^{\!\vee}).\end{cases}
\leqno(2.1.7)$$
Combined with (\hl{2.1.4}{2.1.4}) and Lemma\,\,\hl{L1.2a}{1.2a}, this implies for $\al\in(-1,0\1]$ that
\htt{2.1.8}{}
$$q(\M_{\ee(\al)},F)\les q_{\al}:=\begin{cases}[1\mi\alt_f]&(\{\alt_f\}^{\!\vee}\less\al\,{<}\,0),\\ [-\alt_f]&(\al\,{<}\,\{\alt_f\}^{\!\vee}\,\,\h{or}\,\,\,\al\eq0).\end{cases}
\leqno(2.1.8)$$
\sk
Applying Corollaries~1.1a--b, we can now show that
\htt{2.1.9}{}
$$\Gr_{\Vt}^{\al}\OO_X=\msum_{i=1}^{d_X}\,(\Gr_{\Vt}f_i)\1\Gr_{\Vt}^{\al-1}\OO_X\q\h{if}\,\,\,\al\,{>}\,d_X\mi\alt_f.
\leqno(2.1.9)$$
(Here the situation is rather complicated, since the filtration $\Vt$ on $\OO_X$ is indexed by $\Q$, while the filtration $V$ on $\D_X$ by $\Z$. The {\it strict inequality\1} $\al>d_X\mi \alt_f$ in (\hl{2.1.9}{2.1.9}) is rather essential.) We have to determine which $\al\in\Q$ corresponds to $p$ by (\hl{1.1.11}{1.1.11}) in Corollary\,\,\hl{C1.1a}{1.1a} (that is, $[d_X{-}\al]\eq p$, $\ee(-\al)\eq\ee(\al')$) in the case $-p=q_{\al'}$ so that we have by (\hl{2.1.8}{2.1.8})
\htt{2.1.10}{}
$$\aligned&{-}p=[1\mi\alt_f],\q\{\alt_f\}^{\!\vee}\less\al'\,{<}\,0,\\
\h{or}\q&{-}p=[-\alt_f],\q-1\,{<}\,\al'\,{<}\,\{\alt_f\}^{\!\vee}\,\,\h{or}\,\,\,\al'\eq0.
\endaligned
\leqno(2.1.10)$$
\sk
In the first case of (\hl{2.1.10}{2.1.10}), we have $\alt_f\mi p\in(0,1]$ (that is, $\alt_f\in(p,p{+}1]$), and the condition $[d_X{-}\al]\eq p$ is equivalent to that $\al\in(d_X\mi p\mi 1,d_X\mi p]$. We then get that
\htt{2.1.11}{}
$$\al\in(d_X{-}p{-}1,d_X\mi\alt_f].
\leqno(2.1.11)$$
\sk
In the second case, a similar argument deduces that
\htt{2.1.12}{}
$$\al\in(d_X\mi\alt_f\mi 1,d_X{-}p{-}1].
\leqno(2.1.12)$$
These imply the assertion (\hl{2.1.9}{2.1.9}). The inclusion (\hl{5}{5}) then follows from Corollary\,\,\hl{C1.1d}{1.1d}. This finishes the proof of Theorem~\hl{T2}{2}.
\msn
{\bf 2.2.~Isolated singularity case.} If $Z$ has only one singular point $x$, it is well-known that the action of $f$ on $\OO_{X,x}/(\dd f)$ {\it preserves the filtration $V$ up to the shift by\1} 1, that is,
\htt{2.2.1}{}
$$f\1V^{\al}\bl(\OO_{X,x}/(\dd f)\br)\subset V^{\al+1}\bl(\OO_{X,x}/(\dd f)\br)\q(\al\in\Q).
\leqno(2.2.1)$$
Here the filtration $V$ on $\OO_{X,x}/(\dd f)$ is the quotient filtration of the $V$-filtration on the Brieskorn lattice $H''_f$, see \cite{Bri}, \cite{bl}, \cite{SS}, \cite{Va2}, etc. This filtration coincides with the quotient filtration of the microlocal $V$-filtration $\Vt$ on $\OO_{X,x}$, see \cite[Prop.\,1.4]{JKSY1}.
\sk
As in Corollary\,\,\hl{C1.1a}{1.1a}, this theory gives the {\it canonical\1} isomorphisms
\htt{2.2.2}{}
$$\Gr_V^{\al}\bl(\OO_{X,x}/(\dd f)\br)=\Gr_F^pH^{d_X-1}(F_{\!f},\C)_{\la}\q\bl([d_X\mi\al]=p,\,\,\la=e^{-2\pi i\al}\br),
\leqno(2.2.2)$$
choosing a trivialization of $\omega_X$, where the action of $\Gr_Vf$ on the left-hand side can be identified up to sign with the action of $\Gr_F\1N/2\pi i$ on the right-hand side (implying (\hl{7}{7})), see also Remark\,\,\hl{R1.1a}{1.1a}. Here $F_{\!f}$ denotes the Milnor fiber, and $_{\la}$ means the $\la$-eigenspace for the monodromy (which is the {\it inverse\1} of the Milnor monodromy, see for instance \cite{DiSa}).
\sk
Varchenko's asymptotic Hodge filtration uses essentially the isomorphisms
\htt{2.2.3}{}
$$t^k:\Gr_V^{\al-k}(H''_f[t^{-1}])\simto\Gr_V^{\al}(H''_f[t^{-1}])\q(\al{-}k\,{\in}\,[0,1)),
\leqno(2.2.3)$$
with $H''_f[t^{-1}]$ the localization of the Brieskorn lattice $H''_f$ by $t$, since he studies the coefficient of the initial term of asymptotic expansion of each period integral, see \cite{Va2}. This filtration coincides with the Hodge filtration $F$ (defined by using $\dd_t^{-1}$ on the Gauss-Manin system $H''_f[\dd_t^{-1}]$) after taking the graded pieces of the weight filtration $W$, since the latter is given by the monodromy filtration up to shift, and the action of $N$ vanishes on each graded piece of $W$, see also a remark just after (\hl{1.1.12}{1.1.12}) in Remark~\hl{R1.1a}{1.1a}. 
\par\htt{R2.2a}{}\msn
{\bf Remark\,\,2.2a.} It has been remarked by M.~Musta\c t\u a that the assertion (\hl{4}{4}) in the {\it algebraic case\1} can be reduced to the isolated singularity case using a sequence $h_j\,{:=}\,f\pl c_j\1g_j$ $(j\,{\gg}\,0$). Here the $g_j$ are homogeneous polynomials of degree $j$ having an isolated singularity at $x\in X$ with respect to some local coordinates of $(X,x)$, and $c_j\in\C$. Note that $\{f\pl c_j\1g_j\eq0\}\subset X$ has an isolated singularity at $x$ for a sufficiently general $c_j$. (To see this, consider the intersection
$$\mcap_{i=1}^{d_X}\,\bl\{(1{-}u)\dd f/\dd x_i\pl u\1\dd g_j/\dd x_i\eq0\br\}\,\,\subset\,\,X{\times}\C,$$
with $u$ the coordinate of $\C$. If one takes its intersection with a divisor defined by $u\eq 1$, its dimension decreases at most by 1 on a neighborhood of the divisor.)
\sk
Since $h_j\mi f\in\mm_{X,x}^j$ with $\mm_{X,x}\subset\OO_{X,x}$ the maximal ideal, one can then apply \cite[Prop.\,6.6(3)]{MuPo} (using \cite[Thm\,3.2 or Cor.\,3.3 rather than Thm.\,2.2]{MSS}), which shows that
\htt{2.2.4}{}
$$\lim{}_{j\to\infty}\,\alt_{h_j,x}=\alt_{f,x}.
\leqno(2.2.4)$$
Since $h_j^k\mi f^k\in\mm_{X,x}^j$ and $(\dd h_j)_x\subset(\dd f)_x\pl\mm_{X,x}^{j-1}$, the assertion (\hl{4}{4}) is thus reduced to the isolated singularity case using Krull's intersection theorem
\htt{2.2.5}{}
$$\mcap_{j\gg 0}\,\bl((\dd f)_x\pl\mm_{X,x}^{j-1}\br)=(\dd f)_x.
\leqno(2.2.5)$$
(It does not seem completely trivial to extend the above argument to the analytic case, since certain functors between the derived categories of mixed Hodge modules as in the base change theorem \cite[4.4.3]{mhm} seem to be used implicitly in the proof of \cite[Prop.\,6.6(3)]{MuPo}.)
\par\htt{R2.2b}{}\msn
{\bf Remark\,\,2.2b.} It seems rather difficult to follow some argument in \cite{Sch2} (see also the abstract of \cite{vSt}); for instance, Lemma\,4.2 does not seem to hold if the increasing filtration $\Hst''_k$ is defined by $\dd_t^{k-1}\Hs''\cap\Hst''$ for $k\gess 1$ as in (3.6), since the latter is a {\it lattice\1} of the Gauss-Manin system (using $\dd_t^{-1}\Hs''\eq\Hs'\,{\subset}\,\Hs''$) and its restriction to $T'\,{:=}\,T\,{\setminus}\,\{0\}$ coincides with that of $\Hs''$, and is independent of $k\gess 1$ (shrinking $T$ if necessary), although a similar property does not hold for the Hodge filtration $\Ls^{n+1-k}$ on a locally free sheaf $\Ls$ associated with a compactification of $f$. If there is a splitting $\,s:\Hst''\into\Ls\,$ of $\,i^*:\Ls\onto\Hst''\,$ as in Lemma\,4.2\,(ii), then we would have
\htt{2.2.6}{}
$$s(\Hst''_k)|_{T'}=s(\Hst''_1)|_{T'}\subset\Ls^n|_{T'}=\F^n\q(k\ges 1),
\leqno(2.2.6)$$
since $s$ is a morphism of sheaves. This would imply by Lemma\,4.2\,(i) (using $i^*\ssc s\eq{\rm id}$) that
\htt{2.2.7}{}
$$s(\Hst''_k)\subset j_*\F^n\cap\Ls=\Ls^n,\q\Hst''_k\subset i^*\Ls^n=\Hst''_1=\Hs'',
\leqno(2.2.7)$$
contradicting the above definition of $\Hst''_k\,$ if $\,\Hst''\,{\ne}\,\Hs''$. (Something seems to be misstated.)
\sk
Historically Scherk's conjecture \cite{Sch1} (motivated by a question of Brieskorn, see \cite{vSt}) seems to have stimulated Varchenko's theory of asymptotic mixed Hodge structure (see \cite{Va1}, \cite{Va2}), and the latter inspired some part of the theory of mixed Hodge modules; for instance, the $V$-filtration indexed by $\Q$ and the theory of strict bifiltered complexes for the two filtrations $F,V$, see \cite[Prop.\,2.1]{hf2}, \cite[1.7]{bl}, \cite[Prop.\,3.2]{exp}, \cite[Lem.\,3.3.3--5]{mhp}, etc. The first versions of \cite{exp}, \cite{bl} were written in 1983 during a stay in Grenoble, where Scherk and Steenbrink also visited for some period.
\par\htt{R2.2c}{}\msn
{\bf Remark\,\,2.2c.} It is informed that some argument in \cite{EL} works only in the case $\alt_f\less 1$.
\par\htt{R2.2d}{}\msn
{\bf Remark\,\,2.2d.} The inequality (\hl{7}{7}) does not necessarily hold in the non-isolated singularity case even if we use the vanishing cohomology groups $\Ht^{\ssb}(F_{\!f,x},\Q)$ instead of the vanishing cycle complex $\varphi_f\Q_X$ for the definition of ${\rm NO}(f)$; for instance, consider $f\,{:=}\,gt$ on $X\,{:=}\,Y{\times}\C^*$ with $t$ the coordinate of $\C$ and $g$ any holomorphic function on a complex manifold $Y$ with an isolated singularity satisfying ${\rm NO}(g)\ne 1$.
Note that taking the cokernel of the actions of $\Gr^F_1\D_X$ on $\Gr^F_{\ssb}\M_{\la}$ does not necessarily commute with taking the graded pieces of the weight filtration $W$. In the isolated singularity case, this cokernel is already taken. In the above example, we have $\dd f/\dd t=f/t$, and $t$ is invertible on $X$.

\end{document}